\documentclass[11pt]{amsart}

\usepackage{setspace}
\usepackage{tikz}
\usepackage{booktabs}
\usepackage{hyperref}
\usepackage{longtable}

\usetikzlibrary{trees}
\usepackage{mathrsfs}
\usepackage[english]{babel}
\usepackage[latin1]{inputenc}
\usepackage{amsmath,amssymb,amsthm,epsfig}
\usepackage{color}

\usepackage{latexsym}
\usepackage{amsthm}
\usepackage{tikz}

\usepackage{textcomp}
\usetikzlibrary{trees}

\providecommand{\U}[1]{\protect\rule{.1in}{.1in}}
\providecommand{\U}[1]{\protect\rule{.1in}{.1in}}
\providecommand{\U}[1]{\protect\rule{.1in}{.1in}}
\providecommand{\U}[1]{\protect\rule{.1in}{.1in}}
\providecommand{\U}[1]{\protect\rule{.1in}{.1in}}

\newtheorem{Theorem}{Theorem}[section]

\newtheorem*{conjecturen}{Conjecture}

\newtheorem{Corollary}[Theorem]{Corollary}

\theoremstyle{definition}

\newtheorem{Definition}[Theorem]{Definition}

\numberwithin{equation}{section}

\newcommand{\Res}{\operatorname{Res}}

\newcommand{\supp}{\operatorname{Supp}}

\newcommand{\codim}{\operatorname{codim}}

\newcommand{\Var}{\operatorname{Var}}

\usepackage{times, amscd}

\newcommand{\sO}{\mathscr{O}}

\newcommand{\OO}{\mathcal{O}}

\def \fol {{\mathscr{F}}}
\def \sing {{\rm Sing}}
\def \singf {{\rm Sing}(\mathscr{F})}

\newcommand{\A}{\mathscr{O}}

\newcommand{\CC}{{\mathbb{C}}}

\newcommand{\PP}{{\mathbb{P}}}

\newcommand{\F}{\mathscr{F}}
\newcommand{\M}{\mathcal{M}}

\begin{document}
\title{Brunella--Khanedani--Suwa variational residues for invariant currents}
\author[Mauricio Corr\^ea]{Mauricio Corr\^ea}
\address{\sc Mauricio Corr\^ea\\
UFMG\\
Avenida Ant\^onio Carlos, 6627\\
30161-970 Belo Horizonte\\ Brazil}
\email{mauriciojr@ufmg.br}
\author[Arturo Fern\'andez-P\'erez]{Arturo Fern\'andez-P\'erez}
\address{\sc Arturo Fern\'andez-P\'erez\\
UFMG\\
Avenida Ant\^onio Carlos, 6627\\
30161-970 Belo Horizonte\\ Brazil}
\email{fernandez@ufmg.br}
\author[Marcio G. Soares]{Marcio G. Soares}
\address{\sc Marcio G. Soares\\
UFMG\\
Avenida Ant\^onio Carlos, 6627\\
31270-901 Belo Horizonte\\ Brazil}
\email{msoares@impa.br}

\dedicatory{\it To  Israel Vainsencher on the occasion of his 70th birthday}

\begin{abstract}
In this work we prove a Brunella--Khanedani--Suwa variational  type residue  theorem for  currents invariant by holomorphic foliations. 
As a consequence, we give conditions for the leaves of a singular holomorphic foliation to accumulate in the intersection of the singular set of the foliation with the support of an invariant current. 
\end{abstract}
\maketitle

\setcounter{tocdepth}{1}
\tableofcontents

\section{Introduction}

In \cite{KS} B. Khanedani and T. Suwa 
introduced    an index for singular holomorphic foliations on complex compact  surfaces called the \textit{Variational index}.
In \cite{LS} D. Lehmann and T. Suwa generalized the variational index for higher dimensional holomorphic foliations. In particular, they showed  that if $V$ is a complex subvariety  invariant by a  holomorphic foliation $\F$ of  dimension $k\geq 1$ on a $n$-dimensional complex compact manifold $X$, then 
$$
 c_1^{n-k}(\det(N^*\F)) \cdot  [V]  = (-1)^{n-k}\sum_{\lambda} \Res_{ c_1^{n-k}}(\F; S_\lambda), 
$$
where  $S_\lambda$ is  a connected component of  $S(\F,V):=(\singf\cap V) \cup \sing(V) $ (here $\singf$ and $\sing(V)$ denote the singular sets of $\F$ and $V$, respectively),  $ [V]$ is the  integration current  of $V$ and $N^*\F$ is the conormal sheaf of $\F$. In the case $X$ is a complex surface  and  $S(\F,V)$ is an isolated set, then for each $p\in S(\F,V)$,
$$-\Res_{ c_1}(\F; p)=\Var(\F,V,p),$$ 
where $\Var(\F,V,p)$ denotes the Variational index of $\F$ along $V$ at $p$, as defined by Khanedani-Suwa in \cite{KS}. 
\par M.  Brunella in \cite{Br4} studied the Khanedani-Suwa variational index and its relations with the GSV and the Camacho-Sad indices.  See also \cite[II,  Proposition 1.2.1]{Mc}.

M. McQuillan, in his proof of  the Green-Griffiths  conjecture (for a projective surface $X$ with $c_1^2(X)> c_2(X)$ ), \cite{Mc}, showed   that 
 if  $ X$ is  a complex surface of general type and $\F$ is a holomorphic foliation on $X$, then $\F$ has no entire leaf which is Zariski dense. See \cite{Dem, PS, GPR,Den} for more details about the Green-Griffiths conjecture and generalizations.  M. Brunella in \cite{Br2} provided an alternative  proof of  McQuillan's result by showing that  if $[T_f]$   is the  Ahlfors current associated to a Zariski dense entire curve $f :\mathbb{C} \to  X$ which is tangent to $\F$, then 
$$
c_1(N^*\F) \cdot  [T_f]  = \sum_{p\in \singf \cap \supp(T_f)}  \frac{1}{2\pi i}[T_f] (\chi_{U_p}d(\phi_p \beta_p)) \leq 0,
$$
where $\chi_{U_p}$ denotes the characteristc function of a neighborhood $U_p$ of $p\in \singf \cap \supp(T_f)$, see section \ref{secao-var} for more details.
\par To continue we consider $\F$ a singular holomorphic foliation of dimension $k\geq 1$ on a complex compact manifold $X$ of dimension at least two. We recall that a positive closed current $T$ in $X$ is \emph{invariant} by $\F$ if $T_{|\F} \equiv 0$, that is, $T(\eta)=0$ for every test form $\eta$ vanishing along the leaves of $\F$, so that $T(\eta)$ depends only on the restriction of $\eta$ to the leaves.  
\par In \cite{Br3} M. Brunella   proved  a more general variational index type Theorem for positive closed currents of bidimension $(1,1)$ invariant by  one-dimensional holomorphic foliations, with isolated singularities, on complex compact manifolds. More precisely,  he showed that  if $T$ is an invariant positive closed current of  bidimension $(1,1)$, then 
$$
c_1(\det(N^*\F)) \cdot  [T] = \sum_{p\in \singf \cap \supp(T)}   \frac{1}{2\pi i}[T] (\chi_{U_p}d(\phi_p \beta_p)). 
$$
Compare this formula with the so called  \textit{asymptotic Chern Class} of  a foliation on a complex surface introduced in \cite{CLS}. 
Moreover,  M.  Brunella showed, in the same work,   that a generic one-dimensional  holomorphic foliation on complex projective spaces has no invariant measure. 
In \cite[Corollary 1.2]{K}  L. Kaufmann showed that  there is no diffuse foliated cycle directed by  embedded Lipschitz laminations of dimension $k\geq n/2$ on $\PP^n$. 
\par We denote the class of a closed  current  $T$ of bidimension $(p,p)$ in the  cohomology group $H^{n-p,n-p}(X)$  by $[T]$. In order to provide a generalization of the above results, we define the residue of $\F$ relative to $T$ along a connected component of the singular set of $\F$, (see for instance Def. \ref{res} in Sect. \ref{secao-var}.). In this work we prove the following result. 
\begin{Theorem} \label{teorema}
Let $\F$ be a holomorphic foliation, of dimension $k\geq 1$, on a complex compact manifold $X$ with $\dim(\singf)\leq k-1$. Write $\singf = \displaystyle\bigcup_{\lambda} Z_\lambda$, a decomposition  into connected components and let $U_{\lambda}$ be a regular neighborhood of $Z_\lambda$. For $p \geq k$, if $T$ is a positive closed current of bidimension $(p,p)$ invariant by $\F$, then
$$
c_1^{p-k+1}(\det (N^*{\F}))\cdot [T]=\sum_{Z_\lambda\subset \supp(T)\cap \singf} \Res(\F, T, Z_\lambda).
$$
\end{Theorem}

\par A compact non-empty subset
$\M \subset X$ is said to be a \textit{ minimal set} for $\fol$ if
the following properties are satisfied
 \begin{enumerate}
\item[(i)] $\M$ is invariant by $\fol$;
\item[(ii)] $\M \cap \sing(\fol)= \emptyset$;
 \item[(iii)] $\M$ is minimal with respect to these properties.
\end{enumerate}
\par The problem of existence of minimal sets for codimension one holomorphic foliations on $\mathbb{P}^n$ was considered by Camacho--Lins Neto--Sad in \cite{CLS}.
To our knowledge, this problem remains open for $n=2$.  If $\fol$ is
a codimension one holomorphic foliation on $\mathbb{P}^n$, with $n\geq 3$,
Lins Neto \cite{Ln} proved that $\fol$ has no minimal sets.

M. Brunella posed in \cite{Br1} the following question:
\begin{conjecturen}\label{conje}
Let $X$ be a compact connected complex manifold of dimension $n\geq 3$, and let $\fol$
be a codimension one holomorphic foliation on $X$  such that 
  $N\fol$ is ample. Then every leaf of $\fol$
accumulates to $\singf$.

\end{conjecturen}
\par In \cite{perrone}, Brunella--Perrone
proved the above Conjecture for codimension-one holomorphic
foliations on a projective manifold with cyclic Picard group.  In \cite{C-P}  the natural conjecture has been  stated:

\begin{conjecturen}[Generalized Brunella's conjecture]\label{conje1}
Let $X$ be a compact connected complex manifold of dimension $n\geq
3$, and let $\fol$ be a holomorphic foliation of codimension $r<n$ on $X$   such that 
  $\det(N\fol)$ is ample. Then every
leaf of $\fol$ accumulates to $\singf$, provided $n\geq 2r+1$.
\end{conjecturen}

The main result  in  \cite{C-P}  suggests that the property of accumulation of the leaves of a foliation $\F$  to the  its  singular set 
 (\textit{or nonexistence of minimal sets of $\fol$}  ) depends on  the existence of strongly $q$-convex spaces  which contains the singularities of $\F$.
In  \cite{CLS} it was proved that  there is no invariant measure  with support on a nontrivial minimal set of  a foliation on $\mathbb{P}^2$. We observe that in  $\mathbb{P}^n$ we have that  $\det(N_{\F})$ is ample  for  every foliation $\F$. The following  Corollary \ref{inter-sing} generalize the result in \cite[Theorem 2]{CLS}.

\begin{Corollary}\label{inter-sing}
Let $\F$ be a holomorphic foliation, of dimension $k\geq 1$, on      a projective manifold $X$  such that  $\dim(\singf)\leq k-1$ and  $\det(N\fol)$  ample. 
Suppose that $h^{n-p,n-p} (X ) = 1$,  for some  $p \geq k$.  If  $T$ is a 
  positive closed current of bidimension  $(p,p)$ invariant by $\F$,  then  $\supp(T)\cap \singf\neq \emptyset $. In particular,   there is no invariant  positive closed current of  bidimension  $(p,p)$  with support on a nontrivial minimal set of $\F$. 
\end{Corollary}

Compare Corollary \ref{inter-sing} with  \cite[Corollary 5.5]{K}.
Since $h^{n-p,n-p} (\PP^n ) = 1$, this result holds for foliations on $\PP^n $, in particular   if $V\subset \mathbb{P}^n$ is an $\F$-invariant complex subvariety, then  $V\cap \singf\neq \emptyset.$ This is the   Esteves--Kleiman's result  \cite[12, Proposition 3.4, pp. 12]{EK}.

We can also apply  Theorem  \ref{teorema}  to the  Ahlfors currents   associated to $f: \mathbb{C}^k\to X$   a holomorphic map of generic maximal rank which
is a leaf of the foliation $\F$.  Fix   a K\"ahler form $\omega$ on $X$. On $\mathbb{C}^k$ we take the homogeneous metric form
$$\omega_0:=dd^c \ln |z|^2,$$ and denote by 
$$\sigma=d^c \ln |z|^2 \wedge \omega_0^{k-1}$$
the Poincar\'e form.  Consider $\eta \in A^{1,1}(X)$ and for any $r>0$  define
$$T_{f,r}(\eta)=\int_0^r \frac{dt}{t} \int_{B_t} f^*\eta \wedge \omega_0^{k-1},$$
where $B_t \subset \mathbb{C}^k$ is the ball of radius $t$. Then we consider the positive currents $\Phi_r \in A^{1,1}(X)'$ defined by
$$\Phi_r(\eta):=\frac{T_{f,r}(\eta)}{T_{f,r}(\omega)}.$$
This gives a family of positive currents of bounded mass from which we can extract a subsequence $\Phi_{r_n}$ which converges to a current  $ [T_f] \in A^{1,1}(X)'$ called the 
 Ahlfors current of $f$,  see \cite[Claim 2.1]{GPR}. 
 
 This construction has been  generalized in \cite{BS}  by Burns--Sibony and \cite{DeT}  by De Th\'elin.  In order to  associate to  $f: \mathbb{C}^k\to X$
a positive closed current of any bidimension $(s, s)$, $1 \leq s \leq k $ (also called Ahlfors currents)  it is necessary to impose some technical conditions.

We obtain  another consequence of Theorem \ref{teorema} as follows:
\begin{Corollary}\label{inter-Ahlfors}
Let $\F$ be a holomorphic foliation, of dimension $k\geq 1$, on      a projective manifold $X$  such that  $\dim(\singf)\leq k-1$ and  $\det(N\fol)$  ample. 
Let $f: \mathbb{C}^k\to  X$ be a holomorphic map of generic maximal rank which
is a leaf of the foliation.
Suppose that $h^{n-p,n-p} (X ) = 1$,  for some  $p \geq k$,  and that  there exists an  Ahlfors current of bidimension $(p,p)$   associated  to $f$. 
Then $\overline{f( \mathbb{C}^k)}\cap \singf\neq \emptyset $. 
\end{Corollary}

  \subsection*{Acknowledgments}
MC was partially supported  by CNPq, CAPES and  FAPEMIG. A. F-P was partially supported by CNPq grant number 427388/2016-3.

\section{Singular holomorphic foliations }

Let  $X$ be a connected compact complex manifold of dimension $n$. To define a (singular) holomorphic foliation $\F$ on $X$ we adopt the following point of view: such a $\F$ is determined by a coherent subsheaf $N^\ast \F$ of the cotangent sheaf $T^\ast X = \Omega^1_X$ of $X$ which satisfies
\begin{itemize}
\item[1)] integrability: $d N^\ast \F \subset N^\ast \F \wedge \Omega^1_X$ and
\item[2)] $\Omega^1_X / N^\ast \F$ is torsion free.
\end{itemize}
The generic rank of $N^\ast \F$ is the codimension of $\F$, the dual $(N^\ast \F)^\ast = N\F$ is the normal sheaf of the foliation and the singular locus of $\F$ is
\begin{equation}\label{sing}
\singf= \{p \in X : (\Omega^1_X / N^\ast \F)_p \,\mathrm{is\, not\, a\, free} \,\A_p -\mathrm{module}\}.
\end{equation}
Condition 2 above implies codim$(\singf) \geq 2$. 

Remark that, on $X \setminus \singf$, we have an exact sequence of holomorphic vector bundles
$$
0 \longrightarrow N^\ast \F \longrightarrow \Omega^1_X \longrightarrow T^\ast \F \longrightarrow 0
$$
and, dualizing
$$
0 \longrightarrow T\F \longrightarrow TX \longrightarrow N\F \longrightarrow 0,
$$
where $T\F$ is called the tangent bundle of $\F$, of  dimension $k=(n - \codim( \F))$. Also, since the singular set has codimension greater than $1$ we have the adjunction formula
$$
KX = K\F \otimes \det (N^\ast{\F}),
$$
where $K\F=\det(T\F)^*$ denotes the canonical bundle of $\F$.

If $\F$ has codimension $(n-k)$ then, by taking the $(n-k)$-th wedge product of the inclusion 
$$N^\ast \F \longrightarrow \Omega^1_X,$$ we get a $(n-k)$-form $\omega$ with coefficients in the line bundle $(\bigwedge^{n-k} N\F)^{\ast \ast} = \det (N\F)$.

\subsection{Holomorphic foliations on complex  projective spaces}

Let $\omega \in \mathrm H^0(\PP^n,\Omega_{\PP^n}^{n-k}(m))$ be the twisted $(n-k)$-form induced by a holomorphic foliation $\F$  of dimension  $k$ on $\PP^n$.

Take a    generic non-invariant  linearly embedded subspace $i:L\simeq \PP^{n-k} \hookrightarrow \PP^n$.  We have an induced  non-trivial section  
$i^*\omega \in \mathrm H^0(L ,\Omega_{L}^{n-k}(m)) \simeq \mathrm H^0( \PP^{n-k} , \OO_{\PP^{n-k}}(k-n-1+m)),$ 
since $\Omega_{\PP^{n-k}}^{n-k}=\OO_{\PP^{n-k}}(k-n-1)$ .  The {\it degree} of $\F$ is defined by  $$\deg(\F):=\deg(Z(i^*\omega))=k-n-1+m.$$
In particular,  $\omega \in \mathrm H^0(\PP^n,\Omega_{\PP^n}^k(\deg(\F)+n-k+1))$.  That is,  $\det(N\F)=\OO_{\PP^n}(\deg(\F)+n-k+1)$ is ample. 

A  holomorphic foliation, of degree $d$, can be induced  by  a polynomial $(n-k)$-form on $\CC^{n+1}$ with homogeneous coefficients of degree $d+1$, see for instance \cite{CMM}.
\section{Variational residue and proof of Theorem \ref{teorema}}\label{secao-var}
Hence, a holomorphic foliation of dimension $k$ is given by a family $(\{V_\mu\}, \{\omega_\mu\})_{\mu \in \Lambda}$, where $\mathcal{V} =\{V_\mu\}_{\mu \in \Lambda}$ is an open cover of $X$ by Stein open sets, $\omega_\mu$ is an integrable holomorphic $(n-k)$-form defined in $V_\mu$ and locally decomposable in $V_\mu \setminus \sing({\F})$. This means that, for each $p \in V_\mu$, there is an open neighborhood $V_p \subset V_\mu$ of $p$ such that 
$$\omega_{\mu|V_{p}} = \omega_1^\mu \wedge \dots \wedge \omega_{n-k}^\mu,$$ where $\omega_j^\mu$ is a holomorphic 1-form and $d \omega_j^\mu \wedge \omega_\mu =0$ for $1 \leq j \leq n-k$.

The integrability condition tells us that, in $V_\mu \setminus \sing({\F})$,  there is a $C^\infty$ $1$-form $\alpha_\mu$ satisfying:
\\
\\
\noindent(i) $d\omega_\mu = \alpha_\mu \wedge \omega_\mu$, for all $\mu \in \Lambda$. $\alpha_\mu$ is not unique, but its restriction to the leaves of $\F$ is, provided $\omega_\mu$ is fixed.
\\
\\
\noindent(ii) $\alpha_\mu$ is of type $(1,0)$ since $\omega_\mu$ is holomorphic and ${\alpha_\mu}_{|\F}$ is holomorphic. This last fact follows from: if we assume that around a regular point the foliation $\F$ is generated by $\partial/\partial z_i$, $i=1, \dots, k$, then $\iota_{\partial/\partial z_i}(d\omega_\mu) = \left(\iota_{\partial/\partial z_i}\alpha_\mu\right) \omega_\mu$. In particular, if $k=1$ then $\alpha_{\mu|\F}$ is closed and $d\alpha_{\mu|\F}=0$.
\\
\\
\noindent(iii) In the overlapping $V_{\mu\nu}$ we have $\omega_\mu=f_{\mu\nu}\omega_\nu$, with $f_{\mu\nu} \in {\sO}^{\ast}(V_{\mu\nu})$ and the cocycle $\{f_{\mu\nu}\}_{\mu, \nu \in \Lambda}$ determines the line bundle $\det (N{\F})$. Hence 
\begin{equation}\label{nor1}
\left(\alpha_\mu - \alpha_\nu - \dfrac{df_{\mu\nu}}{f_{\mu\nu}}\right) \wedge \omega_\mu= 0.
\end{equation}
This shows that $\alpha_\mu - \alpha_\nu - \dfrac{df_{\mu\nu}}{f_{\mu\nu}}$ is a $C^\infty$ local section of the conormal bundle $ N^\ast \F$ of the regular foliation ${\F}_{| X \setminus \sing({\F})}$. Since the sheaf of smooth sections of $N^\ast \F$ is acyclic, we have that there exist $C^\infty$ 1-forms $\gamma_\mu$ in $V_\mu$ satisfying: $\gamma_\mu$ is a local section of $N^\ast \F$ and
$$\alpha_\mu - \alpha_\nu - \dfrac{df_{\mu\nu}}{f_{\mu\nu}} = \gamma_\mu - \gamma_\nu,$$ 
so that 
$$\alpha_\mu - \gamma_\mu = \alpha_\nu - \gamma_\nu + \dfrac{df_{\mu\nu}}{f_{\mu\nu}}.$$
Call $\beta_\mu =\alpha_\mu - \gamma_\mu,$
 hence
\begin{equation}\label{nor2} 
\beta_\mu =\beta_\nu + \dfrac{df_{\mu\nu}}{f_{\mu\nu}}, \; d \beta_\mu = d \beta_\nu \; \mathrm{in} \; V_{\mu \nu},\; d\omega_\mu= \beta_\mu \wedge \omega_\mu \; \mathrm{and} \; d \beta_\mu \wedge \omega_\mu =0. 
\end{equation} 
By the second equality in \ref{nor2}, the 2-forms $\{d\beta_\mu\}$ piece together and we have a global $C^\infty$ 2-form on $X \setminus \sing({\F})$ which we denote by $d \beta$. 

We shall briefly digress on the geometric meaning of this smooth 2-form $d\beta$ (see \cite{C-C} 6.2.4): the first equality in \ref{nor2} tells us that the 1-forms $\{ \beta_\mu\}$ behave as connection matrices of $\det (N {\F})$, in $V_\mu$, for some connection. In this case it's natural to consider the basic connections (in the sense of Bott, see \cite{C-L}). 

Fix a $C^{\infty}$ decomposition $$TX_{|X \setminus \sing({\F})} = N{\F} \oplus T{\F},$$
 where $N_{\F}$ and $T_{\F}$ are the normal and tangent bundles, respectively, of the regular foliation ${\F}_{| X \setminus \sing({\F})}$. 

Let $V_\mu$ be the domain of a local trivialization of $N{\F}$ and $\{v_1^\mu, \dots, v_{n-k}^\mu\}$ be a local frame for $N{{\F}_{|V_{\mu}}}$ such that $\omega_\mu (v_1^\mu, \dots, v_{n-k}^\mu) \equiv 1$. For a suitable basic connection $\nabla$ and $\zeta$ any section of $T{{\F}_{|V_\mu}}$, we have that  
$$\beta_\mu(\zeta)=-{\mathrm{tr}} (\theta^\mu)(\zeta)$$
 if, and only if, $d\omega_\mu = \beta_\mu \wedge \omega_\mu$,
where $\theta^\mu$ is the connection matrix in $V_\mu$ of $\nabla$ relative to the frame $\{v_1^\mu, \dots, v_{n-k}^\mu\}$. In particular, the 1-forms $\{\beta_\mu\}$ piece together to give a well defined global form $\beta$ on $X \setminus \sing({\F})$. 
It follows that $d\beta =-{\mathrm{tr}}(K_\nabla) = - c_1(K_\nabla)$ where $K_\nabla=\{K_\nabla^\mu\}_{\mu \in \Lambda}$ is the curvature form of $\nabla$ and the class  $d\beta = -c_1({N}{\F})= c_1( \det N^\ast{\F})$.

\begin{Definition}\label{res}
\em{Let $\F$ be a singular foliation of dimension $k\geq 1$, as above, and consider 
 $$\sing({\F})= \displaystyle\bigcup_{\lambda} Z_\lambda$$
a decomposition of its singular locus into connected components.  For $p \geq k$, suppose $T$ is a positive closed current of bidimension $(p,p)$ which is invariant by $\F$.  The residue of $\F$ relative to $T$ along $Z_\lambda$ is
$$
\Res(\F, T, Z_\lambda) = \dfrac{ T\left(d(\varphi_{_{\lambda}} \beta)^{p-k+1}\wedge  \chi_{_{Z_\lambda}} \upsilon_{_{Z_\lambda}} \right)}{\mathrm{vol}(Z_\lambda )} \cdot [Z_\lambda],
$$
where $\chi_{_{Z_\lambda}} $ denotes the characteristic function, $\upsilon_{_{Z_\lambda}}$ is a volume element of $Z_\lambda$ and $\varphi_\lambda : X \longrightarrow \mathbb{R}$ is a $C^\infty$ function satisfying $\varphi_\lambda^{-1}(0) =Z_\lambda$,  $0 < \varphi_\lambda \leq 1$ in $X \setminus Z_\lambda$ and $\varphi_\lambda =1$ in $X \setminus {{U_{\lambda}}}$.}
\end{Definition}

 Now, we are able to prove Theorem \ref{teorema}:

\noindent\textbf{Theorem}. \emph{Let $\F$ be a holomorphic foliation of dimension $k$ on a complex compact manifold $X$ with $\dim(\singf)\leq k-1$. Write $\singf = \displaystyle\bigcup_{\lambda\in L} Z_\lambda$, a decomposition  into connected components and let $U_{\lambda}$ be a regular neighborhood of $Z_\lambda$. For $p \geq k$, if $T$ is a positive closed current of bidimension $(p,p)$ invariant by $\F$ then,
$$
c_1^{p-k+1}(\det (N^*{\F}))\cdot [T]=\sum_{Z_\lambda\subset \supp(T)\cap \singf} \Res(\F, T, Z_\lambda).
$$
}

\begin{proof}   In order to show geometrically that  $$c_1^{p-k+1}(\det (N^*{\F})) \cdot [T]$$ localizes at $\supp(T)\cap \singf$  we will use the concept of regular neighborhood. 
\par An open set ${U_{\lambda}}$, $Z_\lambda \subset {U_{\lambda}} \subset X$, is a \emph{regular neighborhood} of $Z_\lambda$ provided $\overline{{U_{\lambda}}}$ is a (real) $C^0$ manifold of dimension $n$ with boundary $\partial {U_{\lambda}}$. 
Regular neighborhoods can be obtained in the following way: take a Whitney stratification $\mathcal S$ of $Z_\lambda$ and let $W_\lambda$ be any open set containing $Z_\lambda$. By the proof of Proposition 7.1 of \cite{Ma}, we can construct a family of tubular neighborhoods $\{T_{S, \rho_{S}}\}$, with $|T_{S, \rho_{S}}| \subset W_\lambda$ ($\rho_S$ is the tubular function),  of the strata $S$ of $\mathcal{S}$, satisfying the commutation relations which give control data for $\mathcal{S}$: if $S$ and $S'$ are strata with $S < S'$
then
$$\begin{array}{c}
\pi_S \pi_{S'} (p)= \pi_S(p) \\
\rho_S\pi_{S'} (p)= \rho_S(p).
\end{array}
$$
This allows for the construction of ${U_{\lambda}}$ as a subset of $W_\lambda$ and, by shrinking $W_\lambda$, we may assume ${U_{\lambda}} \cap {U_{{\tilde{\lambda}}}}= \emptyset$ for ${{\lambda}} \neq {{{\tilde{\lambda}}}}$. We call $\{U_\lambda\}_{\lambda \in L}$ a \emph{system of regular neighborhoods} of $Z$. 

Let $\{U_\lambda\}_{\lambda \in L}$be a system of regular neighborhoods of $\singf$. Since
$$d(\varphi_{_{\lambda}} \beta) = d \varphi_\lambda \wedge \beta + \varphi_\lambda d\beta = \varphi_\lambda  d\beta 
$$ 
outside $U_\lambda$, $d \beta_{|\F}=0$ in $X \setminus \singf$ and $T$ is $\F$-invariant, we have 
$$\left(T \wedge  \chi_{_{Z_\lambda}} \upsilon_{_{Z_\lambda}}\right)\left(d(\varphi_{_{\lambda}} \beta)^{p-k+1} \right) =0$$ in $X \setminus U_\lambda$. By reducing the tubular function of $U_\lambda$ we conclude that 
$$
\supp \left(T \wedge  \chi_{_{Z_\lambda}} \upsilon_{_{Z_\lambda}}\right)\left(d(\varphi_{_{\lambda}} \beta)^{p-k+1} \right) \subseteq Z_\lambda
$$
which gives
$$
T \left(d(\varphi_{_{\lambda}} \beta)^{p-k+1} \wedge \chi_{_{Z_\lambda}} \upsilon_{_{Z_\lambda}}\right) = \mu_{_{Z_\lambda}} \left[{Z_\lambda}\right](\upsilon_{_{Z_\lambda}})
$$
and 
$$
\mu_{_{Z_\lambda}} = \Res(\F, T, Z_\lambda).
$$
Since $d(\varphi_{_{\lambda}}\beta)$ represents  $c_1(\det N^\ast{\F})$, we get
$$
c_1^{p-k+1}(\det (N^*{\F}))\cdot [T]=\sum_{Z_\lambda\subset \supp(T)\cap \singf} \Res(\F, T, Z_\lambda).
$$

\end{proof}

\subsection{ Proof of Corollaries \ref{inter-sing} and 
\ref{inter-Ahlfors} }
It is enough to prove the Corollary  \ref{inter-sing}. The result 
 is a straightforward consequence of Theorem \ref{teorema}.  In fact, suppose by contradiction that
 $T$ is a closed positive current of bidimension $(p,p)$ invariant by $\F$ and that $\supp(T)\cap \singf =\emptyset $. Then, it follows from Theorem \ref{teorema} that
 $$
c_1^{p-k+1}(\det (N^*{\F}))\cdot [T]=0.
$$
Since $h^{n-p,n-p} (X ) = 1$ and $\det (N^*{\F})$ is ample, then $[T]=b \cdot c_1^{n-p}(\det (N{\F})) \in H^{n-p,n-p} (X )$, for some $b>0$. 
Therefore, we have 
$$c_1^{p-k+1}(\det (N^*{\F}))\cdot [T]= (-1)^{p-k+1} b\cdot c_1^{n-k+1}(\det (N{\F}))\neq 0.$$
This is a contradiction.

\end{document}